\documentclass{article}
\usepackage{amsmath,amsthm,amssymb}

\newcommand{\kkk}{[\kappa]^{\omega}}

\newcommand{\dom}{\mathrm{dom}}

\newcommand{\force}{\Vdash}
\newcommand{\forces}{\force}

\def\hhh{{[H_\lambda]^\omega}}
\def\kkp{[H_{\kappa^+}]^\omega}

\newtheorem{theorem}{Theorem}[section]
\newtheorem{lemma}[theorem]{Lemma}

\newtheorem{corollary}[theorem]{Corollary}

\theoremstyle{definition}
\newtheorem{definition}[theorem]{Definition}
\newtheorem{example}[theorem]{Example}

\title{On the structure of stationary sets\footnote{2000 AMS subject
classification 03E40, 03E65}} \author{
Qi Feng
\thanks{Partially supported by CNSF Grants 19925102 and 19931020}\\
Institute of Mathematics\\
Chinese Academy of Sciences\\
and\\
Department of Mathematics\\
National University of Singapore
\and
Thomas Jech
\thanks{Supported in part by the GA \v CR grant number 201/02/857}\\
Mathematical Institute ASCR
\and
Jind{\v r}ich Zapletal
\thanks {Partially supported by grants GA \v CR 201-03-0933, NSF DMS-0071437, and visiting appointments at
Chinese Academy of Sciences in Beijing and
National University of Singapore}\\
University of Florida}

\begin{document}

\maketitle
\begin{abstract}
We isolate several classes of stationary sets of
$\kkk$ and investigate implications among them.
Under a large cardinal assumption, we prove a structure theorem for
stationary sets.
\end{abstract}

\section{Introduction}

We investigate stationary sets in the space $\kkk$ of countable subsets
of an uncountable cardinal. We concentrate on the following particular
classes of stationary sets:
\vskip5truept
\hskip4truecm full

\hskip2truecm {$\displaystyle{\nearrow\hskip4truecm \searrow}$}

\centerline{club $\to$ local club $\to$ reflective $\to$
projective stationary}

\hskip2truecm {$\displaystyle{\searrow\hskip4truecm \nearrow}$}

\hskip3.7truecm spanning
$$ \hbox{\rm Fig. 1.1}$$
In the diagram, the $\to$ represents implication. We show among others
that under suitable large cardinal assumption (e.g., under Martin's
Maximum), the diagram collapses to just two classes:
$$\left\{ \aligned
{\hbox{\rm club}}\cr
{\hbox{\rm local club}}\cr
{\hbox{\rm full}}\cr
\endaligned \right\}
\longrightarrow \left\{ \aligned
                               {\hbox{\rm projective stationary}}\cr
                               {\hbox{\rm spanning}}\cr
                               {\hbox{\rm reflective}}\cr
                             \endaligned\right\}$$
$$\hbox{\rm Fig.\ } 1.2$$

\noindent Under the same large cardinal assumption, we prove a structure
theorem for stationary sets: for every stationary set $S$ there exists a
stationary set $A\subset \omega_1$ such that $S$ is spanning above $A$
and nonstationary above $\omega_1-A$.

We also investigate the relation between some of the above properties of
stationary sets on the one hand, and properties of forcing on the other,
in particular the forcing that shoots a continuous $\omega_1$--chain
through a stationary set. We show that the equality of the classes of
projective stationary sets and spanning sets is equivalent to the
equality of the class of stationary--set--preserving forcings and the
class of semiproper forcings.

The work is in a sense a continuation of the previous work
\cite{FJ1} and \cite{FJ2} of the first two authors, and ultimately
of the groundbreaking work of \cite{FMSh} of Foreman, Magidor and
Shelah.

\section{ Definitions}

We work in the spaces $\kkk$ and $[H_\lambda]^\omega$, where
$\kappa$ and $\lambda$ are uncountable cardinals. The concept of a
closed unbounded set and a stationary set has been generalized to
the context of these spaces (cf. \cite{J1}) and the generalization
gained considerable prominence following the work \cite{S} of
Shelah on proper forcing.

The space $\kkk$ is the set of all countable subsets of $\kappa$,
ordered by inclusion; similarly for $\hhh$, where $H_\lambda$ denotes
the set of all sets hereditarily of cardinality less than $\lambda$. A
set $C$ in this space is {\sl closed unbounded (club)} if it is closed
under unions of increasing countable chains, and cofinal in the ordering
by inclusion. A set $S$ is {\sl stationary} if it meets every club set.
We shall (with some exceptions) only consider $\kappa$ and $\lambda$
that are greater than $\omega_1$; note that the set $\omega_1$ is a club
in the space $[\omega_1]^\omega$ (which motivates the generalization).
In order to simplify some statements and some arguments, we shall only
consider those $x\in \kkk$ (those $M\in\hhh$) whose intersection with
$\omega_1$ is a countable ordinal (these objects form a club set); we
denote this countable ordinal by $\delta_x$ or $\delta_M$ respectively:

\vskip10truept \noindent$(2.1) \hskip 1truecm \delta_x = x\cap
\omega_1,\;\; \delta_M = M\cap \omega_1.$ \vskip10truept

The filter generated by the club sets in $\kkk$ is generated by the club
sets of the form

\vskip10truept
\noindent$(2.2) \hskip 1truecm C_F=\{x\;|\;x \hbox{\rm \ is closed under
} F\}$
\vskip10truept

\noindent where $F$ is an operation, $F : \kappa^{<\omega}\to\kappa;$
similarly for $H_\lambda$. In the case of $H_\lambda$, we consider only
those $M\in\hhh$ that are submodels of the model $(H_\lambda,\in,<)$,
where $<$ is some fixed well ordering; in particular, the $M$'s are
closed under the canonical Skolem functions obtained from the well
ordering.

For technical reasons, when considering continuous chains in $\kkk$ or
$\hhh$, we always assume that when $\langle x_\alpha\;|\;\alpha <
\gamma\rangle$ is such a chain then for every $\alpha,\beta<\gamma$,

\vskip10truept
\noindent$(2.3)$ \hskip 1truecm if $\alpha<\beta$ then
$\delta_{x_\alpha}<\delta_{x_\beta}.$
\vskip10truept

\noindent The term $\omega_1$--{\sl chain} or $(\gamma+1)$--{\sl chain}, where
$\gamma<\omega_1$, is an abbreviation for ``a continuous
$\omega_1$--chain that satisfies $(2.3)$.''

We also note that in one instance we consider club (stationary) sets in
the spaces $[\kappa]^{\omega_1}$ (where $\kappa\geq \omega_2$) those are
defined appropriately.

Throughout the paper we employ the operations of {\sl projection} and
{\sl lifting}, that move sets between the spaces $\kkk$ for different
$\kappa$:

If $\kappa_1<\kappa_2$ and if $S$ is a set in $[\kappa_2]^\omega$, then
the projection of $S$ to $\kappa_1$ is the set

\vskip10truept
\noindent$(2.4) \hskip 1truecm \pi(S) = \{x\cap \kappa_1\;|\;x\in S\}.$
\vskip10truept

If $S$ is a set in $[\kappa_1]^\omega$ then the lifting of $S$ to
$\kappa_2$ is the set

\vskip10truept
\noindent$(2.5) \hskip1truecm \hat S = \{x\in
[\kappa_2]^\omega\;|\;x\cap\kappa_1\in S\}.$
\vskip10truept

We recall that the lifting of a club set is a club set and the
projection of a club set contains a club set. Hence, stationarity is
preserved under lifting and projection.

The special case of projection and lifting is when $\kappa = \omega_1$:

\vskip10truept
$\pi(S) = \{\delta_x\;|\;x\in S\},\hskip10truept \hat A = \{x\;|\;\delta_x\in A\}
\hfill {(A\subset \omega_1).}$

\begin{definition}
A set $S\subset\kkk$ is a {\sl local club} if the set
$$\{X\in[\kappa]^{\aleph_1}\;|\;S\cap[X]^\omega \hbox{\rm \ contains a
club in } [X]^\omega\}$$
contains a club in $[\kappa]^{\aleph_1}.$
\end{definition}

\begin{definition}
A set $S\subset\kkk$ is {\sl full} if for every stationary
$A\subset\omega_1$ there exist a stationary $B\subset A$ and a
club $C$ in $\kkk$ such that
$$\{x\in C\;|\;\delta_x\in B\}\subset S.$$
(``$S$ contains a club above densely many stationary $B\subset
\omega_1$.'')
\end{definition}

\begin{definition}
A set $S\subset\kkk$ is {\sl projective stationary} if for every
stationary set $A\subset \omega_1$, the set $\{x\in S\;|\;\delta_x\in
A\}$ is stationary.
\noindent(``$S$ is stationary above every stationary $A\subset
\omega_1$.'')
\end{definition}

\begin{definition}
A set $S\subset\kkk$ is {\sl reflective} if for every club $C$ in
$\kkk$, $S\cap C$ contains an $\omega_1$--chain.
\end{definition}

\begin{definition}
%
If $x$ and $y$ are in $[\kappa]^\omega$, then $y$ is an
$\omega_1$--{\sl extension}
of $x$ if $x\subset y$ and $\delta_x = \delta_y$.
\end{definition}

\begin{definition}
A set $S\subset\kkk$ is {\sl spanning} if for every
$\lambda\geq\kappa$, for every club set $C$
in $[\lambda]^\omega$ there exists a
club $D$
in $[\lambda]^\omega$ such that every $x\in D$ has an
$\omega_1$--extension $y\in C$ such that $y\cap \kappa \in S.$
\end{definition}

Local clubs were defined in \cite{FJ1}. Projective stationary sets were
defined in \cite{FJ2}; so were full sets (without the name). Note that
all five properties defined are invariant under the equivalence mod
club filter. All five properties are also preserved under lifting
and projection. For instance, let $S\subset [\kappa_1]^\omega$ be
reflective and let us show that the lifting $\hat S$ to
$[\kappa_2]^\omega$ is reflective. Let $C$ be a club set in
$[\kappa_2]^\omega$ and let $F :\kappa_2^{<\omega}\to\kappa_2$ be such
that $C_F\subset C$. If we let for every $e\in [\kappa_1]^{<\omega}$,
$$f(e) = \kappa_1\cap cl_F(e),$$
where $cl_F(e)$ is the closure of $e$ under $F$, then $C_f$ is a club in
$[\kappa_1]^\omega$. Also for every $x\in C_f$, if $y$ is the closure of
$x$ under $F$ then $y\cap \kappa_1 = x$. Let $\langle
x_\alpha\;|\;\alpha<\omega_1\rangle$ be an $\omega_1$--chain in $S\cap
C_f$, we then let $y_\alpha$ be the closure of $x_\alpha$ under $F$,
then $\langle y_\alpha\;|\;\alpha<\omega_1\rangle$ is an
$\omega_1$--chain in $\hat S \cap C_F$.
The arguments are simpler for the other four properties as well as for
projection.

It is not difficult to see that all the implications in Fig 1.1 hold.
For instance, to see that every spanning set is projective stationary,
note that the definition of projective stationary can be reformulated as
follows: for every club $C$ in $\kkk$, the projection of $S\cap C$ to
$\omega_1$ contains a club in $\omega_1$. So let $C$ be a club in
$\kkk$.
If $S$ is spanning, then there is a club $D$
in $\kkk$ such that all $x\in D$ have an $\omega_1$--extension in $
S \cap C$. Hence $\pi(D)\subset \pi(S\cap C)$, where $\pi$ denotes
the projection to $\omega_1$.

\section{Local clubs and full sets}

Local clubs form a $\sigma$--complete normal filter that extends the
club filter. Local clubs need not contain a club, but they do under the
large cardinal assumption {\sl Weak Reflection Principle} (WRP).

\begin{definition}\cite{FMSh} Weak Reflection Principle at $\kappa$: for
every stationary set $S\subset \kkk$ there exists a set $X$ of size
$\aleph_1$ such that $\omega_1\subset X$ and $S\cap [X]^\omega$ is
stationary in $[X]^\omega$ ($S$ reflects at $X$).
\end{definition}

It is not hard to show \cite{FJ1} that WRP at $\kappa$ implies a
stronger version, namely that for every stationary set $S\subset\kkk$,
the set of all $X\in [\kappa]^{\omega_1}$ at which $S$ reflects is
stationary in $[\kappa]^{\omega_1}$. In other words, every local club in
$\kkk$ contains a club.

Thus WRP is equivalent to the statement that every local club contains a
club. And clearly, WRP at $\lambda > \kappa$ implies WRP at $\kappa$.
The consistency strength of WRP at $\omega_2$ is exactly that of the
existence of a weakly compact cardinal; the consistency of full WRP is
considerably stronger but not known exactly at this time.

\begin{example}
For every ordinal $\eta$ such that $\omega_1\leq \eta<\omega_2$, let
$C_\eta$ be a club set of $[\eta]^\omega$ of order--type $\omega_1$
(therefore $|C_\eta| = \aleph_1$). Let
$\displaystyle{ S = \bigcup\{C_\eta\;|\;\omega_1\leq \eta<\omega_2\} }.$
Then $S$ is a local club in $[\omega_2]^\omega$ and has cardinality
$\aleph_2$. By a theorem of Baumgartner and Taylor \cite{B}, every club
set
in $[\omega_2]^\omega$ has size $\aleph_2^{\aleph_0}$. Therefore, WRP at
$\omega_2$ implies $2^{\aleph_0}\leq \aleph_2$, a result of Todor\v
cevi\' c \cite {T}.
\end{example}

Let $P$ be a notion of forcing and assume that $|P|\geq\aleph_1$. Let
$\lambda\geq|P|^+$ and consider the model $H_\lambda$ whose language has
predicates for forcing $P$ as well as the forcing relation. Note that
every countable ordinal has a $P$--name in $H_\lambda$.

If $M\in \hhh$, a condition $q$ is {\sl semi--generic} for $M$ if for
every name $\dot\alpha$ for a countable ordinal such that $\dot\alpha\in
M$ there exists some $\beta\in M$ such that $q\forces \dot\alpha=\beta$.

The forcing $P$ is {\sl semiproper} (Shelah \cite{S}) if the set

\vskip10truept \noindent$(3.1)\hskip1truecm \{M\in\hhh\;|\;\forall
p\in M\;\exists q<p\; q \hbox{\rm \ is semigeneric for } M\}$
\vskip10truept

\noindent contains a club in $\hhh$.

In \cite{FJ1}, it is proved that $P$ preserves stationary sets (in
$\omega_1$) if and only if the set $(3.1)$ is a local club. Since
$|H_{|P|^+}| = 2^{|P|}$, we conclude that if $P$ is
stationary--set--preserving, then WRP at $2^{|P|}$ implies that $P$ is
semiproper. Consequently, we have

\begin{theorem}\cite{FMSh}
WRP implies that the class of stationary--set--preserving forcing
notions equals the class of semiproper forcing notions.
\end{theorem}

\begin{example}
Namba forcing \cite{N}. This is a forcing (of cardinality
$2^{\aleph_2}$) that adds a countable cofinal subset of $\omega_2$
without adding new reals (cf. \cite{jech}). It preserves
stationary subsets of $\omega_1$ and by Shelah \cite{S}, it is not
semiproper unless $0^{\#}$ exists.
\end{example}

We use the Namba forcing to get a partial converse of Theorem 3.3: if
stationary--set--preserving equals semiproper, then WRP holds at
$\omega_2$.

\begin{theorem}
If there exists a stationary set $S\subset[\omega_2]^\omega$ that does
not reflect, then the Namba forcing is not semiproper. Hence if every
stationary--set--preserving forcing of size $2^{\aleph_2}$ is
semiproper, then WRP holds at $\omega_2$, and every local club in
$[\omega_2]^\omega$ contains a club.
\end{theorem}

\begin{proof}
Let $S\subset [\omega_2]^\omega$ be nonreflecting stationary set and
assume that the Namba forcing $P$ is semiproper.

Since $S$ does not reflect, there exists for each $\alpha$,
$\omega_1\leq\alpha<\omega_2$, an operation
$F_\alpha:\alpha^{<\omega}\to\alpha$ such that no $x\in S$ is closed
under $F_\alpha$.

Let $\lambda = (2^{\aleph_2})^+$. As the set $(3.1)$ contains a club,
there exists some $M\in \hhh$ such that $M\cap\omega_2\in S$,
$\langle F_\alpha\;|\;\omega_1\leq\alpha<\omega_2\rangle\in M$ and there
exists some $q\in P$ semigeneric for $M$.

Let $G$ be a $P$--generic filter (over $V$) such that $q\in G$. In
$V[G]$, look at $M[G]$, where $M[G] = \{\dot x/G\;|\;\dot x\in
M\}$. Since $G$ produces a countable cofinal subset of
$\omega_2^{V}$, $M[G]\cap \omega_2^{V}$ is cofinal in
$\omega_2^{V}$. Let $\alpha<\omega_2^{V}$ be the least ordinal in
$M[G]$ that is not in $M$. Since $G$ contains a semigeneric
condition for $M$, we have $M[G]\cap\omega_1=M\cap\omega_1$ and so
$\omega_1\leq\alpha<\omega_2^{V}$ and
$M[G]\cap\alpha=M\cap\alpha.$ Since $\alpha\in M[G]$, $F_\alpha\in
M[G]$. Hence $M[G]\cap\alpha$ is closed under $F_\alpha$. It
follows that $x= M\cap\alpha = M[G]\cap\alpha$ belongs to $S$ and
is closed under $F_\alpha$. This is a contradiction.
\end{proof}

Now we turn our attention to full sets. First we reformulate the
definition: $S\subset\kkk$ is full if and only if there exists a
maximal antichain $W$ of stationary subsets of $\omega_1$ such that for
every $A\in W$, there exists a club $C_A$ in $\kkk$ with $\hat A\cap
C_A\subset S$, where $\hat A$ is the lifting of $A$ from $\omega_1$ to
$\kkk$.

We remark that the full sets form a filter, not necessarily
$\sigma$--complete. It is proved in \cite{FJ2} that
$\sigma$--completeness of the filter of full sets is equivalent to the
presaturation of the nonstationary ideal on $\omega_1$. It is also known
that presaturation follows from WRP which shows that WRP is a large
cardinal assumption.

\begin{example}
Let $W$ be a maximal antichain of stationary subsets of $\omega_1$ and
consider the model $\langle H_\lambda,\in,<,\cdots\rangle$, whose
language has a predicate for $W$. Let
$$S_W=\{ M\in \hhh\;|\;(\exists A\in W\cap M)\;\delta_M\in A\}.$$
The clubs $C_A=\{M\in\hhh\;|\;A\in M\}$ for $A\in W$ witness that $S_W$
is full.
\end{example}

We will now show that the sets $S_W$ from Example 3.5 generate the
filter of full sets:

\begin{lemma}
Let $S$ be a full set in $\kkk$. There exists a model $\langle
H_\lambda,\in,<,\cdots\rangle$, where $\lambda = \kappa^+$, and a
maximal antichain $W$ of stationary subsets of $\omega_1$ such
that $S_W\subset \hat S.$
\end{lemma}

\begin{proof}
Let $S$ be full in $\kkk$. By the reformulation of full sets, let $W$ be
 a maximal antichain and for each $A\in W$, let $F_A :
\kappa^{<\omega}\to \kappa$ be an operation such that
$$\{x\in C_{F_A}\;|\;\delta_x\in A\}\subset S.$$
Consider a model $\langle H_\lambda,\in,<,\cdots\rangle$, $\lambda
= \kappa^+$, whose language has a predicate for $W$ as well as for
the function assigning the operation $F_A$ to each $A\in W$. We
claim that for every $M\in S_W$, $M\cap\kappa\in S$. To see this,
let $M\in \hhh$ and let
 $A\in W\cap M$ be such that $\delta_M\in A$. Then $M$ is closed under
$F_A$ and so $M\cap \kappa\in C_{F_A}$ and $\delta_{M\cap\kappa} =
\delta_M\in A$. Hence, $M\cap\kappa\in S$.
\end{proof}

Consequently, the filter of full sets on $\kkk$ is generated by
the projections of the sets $S_W$ on $\hhh$ with $\lambda =
\kappa^+$.

In \cite{FJ2}, it is proved that the statement that every full set
contains a club is equivalent to the saturation of the nonstationary
ideal on $\omega_1$ (and so is the statement that every full set
contains an $\omega_1$--chain). More precisely,

\begin{theorem}\cite{FJ2}
(a) If the nonstationary ideal on $\omega_1$ is saturated then for every
$\kappa\geq\omega_2$, every full set in $\kkk$ contains a club.

(b) If every full set in $[H_{\omega_2}]^\omega$ contains an
$\omega_1$--chain, then the ideal of nonstationary subsets of $\omega_1$
is saturated.
\end{theorem}

Consequently, ``every full set is reflective'' is equivalent to ``every
full set contains a club'' and follows from large cardinal assumptions
(such as MM). The consistency of ``full $=$ club'', being that of the
saturation of $NS_{\omega_1}$, is quite strong. Neither ``local club $=$
club'' nor ``full $=$ club'' implies  the other: WRP has a model in
which $NS_{\omega_1}$ is not saturated, while the saturation of
$NS_{\omega_1}$ is consistent with $2^{\aleph_0}>\aleph_2$ which
contradicts WRP. Both are consequences of MM, which therefore implies
that ``club $=$ local club $=$ full''.

\section{Projective stationary and spanning sets}

In this section, we investigate projective stationary and spanning sets
and particularly a forcing notion associated with such sets. Among
others we show that WRP implies that every projective stationary set is
spanning (and then spanning $=$ projective stationary).

First we prove a theorem (that generalizes Baumgartner and Taylor's
result \cite{B}
on clubs) that shows that the equality does not hold in ZFC. Every
spanning subset of $[\omega_2]^\omega$ has size $\aleph_2^{\aleph_0}$
while Example 3.2 gives a projective stationary (even a local club) set
of $[\omega_2]^\omega$ of size $\aleph_2$. Thus the equality ``spanning
$=$ projective stationary'' implies $2^{\aleph_0}\leq \aleph_2$.

\begin{theorem}
Every spanning set in $[\omega_2]^\omega$ has size
$\aleph_2^{\aleph_0}$.
\end{theorem}

\begin{proof}
Let $S\subset [\omega_2]^\omega$ be spanning. We shall find
$2^{\aleph_0}$ distinct elements of $S$. Let $F : [\omega_2]^2\to
\omega_1$ be such that for each $\eta<\omega_2$, the function $F_\eta$,
defined by $F_\eta(\xi) = F(\{\xi,\eta\})$, is a one--to--one mapping
of $\eta$ to $\omega_1$. As $S$ is spanning, there exists an operation
$G$ on $\omega_2$ such that every
$M\in [\omega_2]^\omega$ closed under $G$ has an $\omega_1$--extension
$N$ that is closed under $F$ and $N\in S$.

We shall find models $M_f$, $f\in 2^\omega$, closed under $G$, and
$\delta<\omega_1$ such that

\vskip10truept
\noindent$(4.1)$
{(a)} $\delta_{M_f}\leq \delta$ for each $f$, and

\hskip8truept{(b)} if $f\not=g$ then there exist $\xi\in M_f$ and
$\eta\in M_g$ such that $F(\xi,\eta)\geq\delta.$
\vskip10truept

Now assume that we have models $M_f$ that satisfy $(4.1)$. If $f\not=g$
and if $x\in [\omega_2]^\omega$ is such that $M_f\cup
M_g\subset x$ and $x$ is closed under $F$, then
$\delta_x>\delta.$ Hence if $N_f$ and $N_g$ are $\omega_1$--extensions
of $M_f$ and $M_g$, respectively, and are closed under $F$, then
$N_f\not=N_g$. Thus we get $\{N_f,\;|\;f\in
2^\omega\}$ such that the $N_f$'s are $2^{\aleph_0}$ elements
of $S$.

Toward the construction of the models $M_f$, let
$c_\alpha\subset\alpha$,
for each $\alpha<\omega_2$ of cofinality $\omega$, be a set of order
type $\omega$ with $\sup c_\alpha = \alpha$ and let $M_\alpha$ be the
closure of $c_\alpha$ under $G$. Let $Z\subset\omega_2$ and
$\delta<\omega_1$ be such that $Z$ is stationary and for each $\alpha\in
Z$, $M_\alpha\subset\alpha$ and
$\delta_{M_\alpha}=\delta.$

We shall find, for each $s\in 2^{<\omega}$ (the set of all finite
$0$--$1$--sequences), a stationary set $Z_s$ and an ordinal
$\xi_s<\omega_2$ such that

\vskip10truept
\noindent$(4.2)$
{(i)} if $s\subset t$, then $Z_t\subset Z_s$,

\hskip8truept{(ii)} $(\forall\;\alpha\in Z_s)\;\xi_s\in c_\alpha$,

\hskip8truept{(iii)} $\xi_{\langle s0\rangle} < \xi_{\langle
s1\rangle}$  and $F(\xi_{\langle s0\rangle}, \xi_{\langle
s1\rangle})\geq\delta.$
\vskip10truept
Once we have the ordinals $\xi_s$, we let, for each
$f\in 2^\omega$, $M_f$ be the closure under $G$ of the set
$\{\xi_{f\restriction n}\;|\;n<\omega\}$. Clearly,
$$M_f=\bigcup_{n=0}^{\infty}M_{f\restriction n},$$ where for each $s\in
2^{<\omega}$, $M_s$ is the closure under $G$ of
$\{\xi_{s\restriction 0},\cdots, \xi_s\}$. Since $M_s\subset
M_\alpha$ for $\alpha\in Z_s$, we have $\delta_{M_f}\leq \delta$
for every $f\in 2^\omega$. The condition $(4.2)(iii)$ guarantees
that the models $M_f$ satisfy $(4.1)$.

The $Z_s$ and $\xi_s$ are constructed by induction on $|s|$. Given
$Z_s$, there are $\aleph_2$ ordinals $\xi$ such that $S_\xi =
\{\alpha\in Z_s\;|\;\xi\in c_\alpha\}$ is stationary. Consider the first
$\omega_1+1$ of these $\xi$'s and let $\eta = \xi_{\langle s1\rangle}$
be the $\omega_1+1$st element, and $Z_{\langle s1\rangle} = S_\eta$.
Then find some $\xi < \eta$ among the first $\omega_1$ elements such
that $F_\eta(\xi)\geq\delta$ and let $\xi_{\langle s0\rangle}$ be such
ordinal $\xi$ and let $Z_{\langle s0\rangle} = S_\xi.$
\end{proof}

In Definition 2.6, we defined spanning sets in $\kkk$ as satisfying a
certain condition at every $\lambda\geq\kappa$. The following lemma
shows that it is enough to consider the condition at $H_{\kappa^+}$.

\begin{lemma}
A set $S\subset\kkk$ is spanning if and only if for every club $C$ in
$[H_{\kappa^+}]^\omega$ there exists a club $D$ in
$[H_{\kappa^+}]^\omega$ such that every $M\in D$ has an
$\omega_1$--extension $N\in C$ such that $N\cap\kappa\in S$.
\end{lemma}

\begin{proof}
It is easy to verify that if the condition $\forall C\;\exists D $
etc. holds at some $\mu>\lambda$ then it holds at $\lambda$. Thus
assume that $\lambda\geq\kappa^+$ and the condition of the lemma
holds and let us prove that for every club $C$ in $\hhh$ there
exists a club $D$ in $\hhh$ such that every $M\in D$ has an
$\omega_1$--extension $N\in C$ such that $N\cap\kappa\in S$.

Let $C$ be a club in $\hhh$ and let $F$ be an operation on $H_\lambda$
such that $C_F\subset C$. Let $C_0$ be a club in
$[H_{\kappa^+}]^\omega$ be such that $\hat C_0\subset C_F$. Let $D_0$
be a club in $[H_{\kappa^+}]^\omega$ such that every $M_0\in D_0$ has an
$\omega_1$--extension $N_0\in C_0$ with $N_0\cap\kappa\in S$. Let
$D=\hat D_0$ be the set of all $M\in \hhh$ such that $M\cap
H_{\kappa^+}\in D_0$. Let $M\in D$ and let $M_0 = M\cap H_{\kappa^+}$.
Then $M_0\in D_0$. Let $N_0\in C_0$ be an $\omega_1$--extension of $M_0$
such that $N_0\cap\kappa\in S$. We let $N$ be the $F$--closure of $M\cup
(N_0\cap\kappa)$ in $H_\lambda$. The model $N$ is in $C_F$. We claim
that $N\cap\kappa=N_0\cap\kappa$. This shall give us that
$N\cap\kappa\in S$ and $N$ is an $\omega_1$--extension of $M$.

Let $\alpha\in N\cap\kappa$. Let $\tau$ be a skolem term in
$(H_\lambda,\in,<,F)$ and let $a\in M$ and $\alpha_0,\cdots,\alpha_n\in
N_0\cap \kappa$ be such that $\alpha =
\tau(a,\alpha_0,\cdots,\alpha_n).$
Define $h:[\kappa]^{n+1}\to\kappa$ by
$$h(\beta_0,\cdots,\beta_n) =
\begin{cases}{\tau(a,\beta_0,\cdots,\beta_n)}&{\hbox{\rm if }
\tau(a,\beta_0,\cdots,\beta_n)<\kappa,}\cr
      { 0}&{\hbox{\rm otherwise.}}
\end{cases}
$$
Then $h\in M$ and hence $h\in M\cap H_{\kappa^+} = M_0\subset N_0.$
Therefore,
$$\alpha =h(\alpha_0,\cdots,\alpha_n)\in N_0\cap\kappa.$$
\end{proof}

\begin{definition}
Let $S\subset \kkk$ be a stationary set. $P_S$ is the forcing notion
that shoots an $\omega_1$--chain through $S$: forcing conditions are
continuous $(\gamma+1)$--chains, $\langle
x_\alpha\;|\;\alpha\leq\gamma\rangle$, $\gamma<\omega_1$, such that
$x_\alpha\in S$ for each $\alpha$, and
$\delta_{x_\alpha}<\delta_{x_\beta}$ when $\alpha<\beta\leq\gamma.$ The
ordering is by extension.
\end{definition}

The forcing $P_S$ does not add new countable sets and so $\omega_1$ is
preserved. The generic $\omega_1$--chain is cofinal in $\kkk$ and so
$\kappa$ is collapsed to $\omega_1$.

The following theorem gives a characterization of projective stationary
sets and spanning sets in terms of the forcing $P_S$:

\begin{theorem}
{ (a)} A set $S\subset \kkk$ is projective stationary if and only
if the forcing $P_S$ preserves stationary subsets of $\omega_1$.

{ (b)} A set $S\subset\kkk$ is spanning if and only if the
forcing $P_S$ is semiproper.
\end{theorem}

\begin{proof}
(a) This equivalence was proved in \cite{FJ2}; we include the proof for
the sake of completeness.

Let $A$ be a stationary subset of $\omega_1$. We will show that $P_S$
preserves $A$ if and only if $\hat A \cap S$ is stationary.

First assume that $\hat A\cap S$ is nonstationary and let
$C\subset\omega_1$ be a club such that for every $x\in S$,
$\delta_x\not\in C\cap A$.

Let $\langle x_\alpha\;|\;\alpha<\omega_1\rangle$ be a generic
$\omega_1$--chain and let $D =
C\cap\{\delta_{x_\alpha}\;|\;\alpha<\omega_1\}.$ Then $D$ is a club in
$V[G]$ disjoint from $A$.

Conversely, assume that $\hat A\cap S$ is stationary. We will show that
$A$ remains stationary in $V[G]$. Let $\dot C$ be a name for a club in
$\omega_1$ and let $p$ be a condition. Let $\lambda$ be sufficiently
large. Since $\hat A\cap \hat S$ is stationary in $\hhh$, there exists a
countable model $M$ containing $\dot C$ and $p$ such that $\delta_M\in
A$ and $M\cap \kappa\in S$. Let $\langle
x_\alpha\;|\;\alpha<\delta_M\rangle$ be an $M$--generic
$\delta_M$--chain extending $p$. By genericity,
$\displaystyle{M\cap\kappa = \bigcup\{x_\alpha\;|\;\alpha<\delta_M\} }.$
Since $M\cap\kappa\in S$, it can be added on top of the chain $\langle
x_\alpha\;|\;\alpha<\delta_M\rangle$ to form a condition $q$. This
condition extends $p$ and forces that $\delta_M$ is a limit point of
$\dot C$, and hence $q$ forces that $\delta_M \in \dot C \cap A.$
Therefore, $A$ is stationary in $V[G]$.

(b) First let $S$ be a spanning set in $\kkk$. Let $\lambda\geq
(2^\kappa)^+$ (note that $|P_S|\leq 2^\kappa$) and let us prove that the
set $(3.1)$ contains a club in $\hhh$.

Let $C$ be the club of all models $N\in \hhh$ that contain $S$, the
forcing $P_S$ and the forcing relation. By definition $2.6$, there exists
a club $D$ in $\hhh$ such that every $M\in D$ has an
$\omega_1$--extension $N\in C$ such that $N\cap\kappa\in S$. We claim
that the set $(3.1)$ contains $D$.

Let $M\in D$ and $p\in M$. Let $N\in C$ be an $\omega_1$--extension of
$M$ such that $N\cap\kappa\in S$. We enumerate
all ordinals in $N\cap \kappa$ and
all names $\dot\alpha \in
N$ for ordinals.
Starting with $p_0=p$, construct a sequence of conditions
$p_0>p_1>\cdots > p_n>\cdots$ such that $p_n\in N$ for each $n$, and
for every $\dot\alpha\in N$ there are
some $p_n$ and $\beta\in N$ such that $p_n\forces \dot\alpha = \beta$,
and that for every $\gamma\in N\cap\kappa$ there is some $p_n = \langle
x_\xi\;|\;\xi\leq\alpha\rangle$ such that $\gamma\in x_\alpha$. The
sequence produces a continuous chain whose limit is the set $N\cap
\kappa$. Since $N\cap \kappa\in S$, it can be put on top of this chain
to form a condition $q<p$ that decides every ordinal name in $N$
as an ordinal in $N$. Now since $N$ is an $\omega_1$--extension
of $M$, they have the same set of countable ordinals and it follows that
$q$ is semigeneric for $M$.

Conversely, assume that $P_S$ is semiproper. Let $\lambda\geq
(2^\kappa)^+$ and let $C$ be a club in $\hhh$. Let $F$ be an operation
on $H_\lambda$ such that $C_F\subset C$. Let $\mu>\lambda$ be such
that $F\in H_\mu$. Since $P_S$ is semiproper, there is a club
$D\subset\hhh$ such that every model in $D$ has the form $M\cap
H_\lambda$, where $F\in M\in [H_\mu]^\omega$, and there is a semigeneric
condition for $M$. We shall prove that every $M\cap H_\lambda\in D$ has
an $\omega_1$--extension $N$ in $C_F$ such that $N\cap \kappa\in S$.

Let $M\cap H_\lambda\in D$ and let $q$ be a semigeneric condition for
$M\in [H_\mu]^\omega$. Let $G$ be a generic filter on $P_S$ over $V$
such that $q\in G$. Working in $V[G]$,
let $M[G]$ be the set of all $\dot a/G$ for $\dot a \in M,$ and
let $N=M[G]\cap (H_\lambda)^V$.
Since $P_S$ does not add new countable sets, $N\in V$.  Since $F\in
M[G]$, $M[G]$ is closed under $F$, and so is $N$. Hence $N\in C_F$.
Since
$q$ is semigeneric for $M$, $M[G]\cap\omega_1 = M\cap\omega_1$, and so
$N$ is an $\omega_1$--extension of $M\cap H_\lambda$. Since
the union of the generic $\omega_1$--chain $\langle
x_\alpha\;|\;\alpha<\omega_1\rangle$ is $\kappa$, we claim that the
union of $\langle
x_\alpha\;|\;\alpha<\delta_M\rangle$ is $M[G]\cap \kappa = N\cap\kappa$.
Granting this claim,
this union is $x_{\delta_M}$ and $\langle
x_\alpha\;|\;\alpha\leq\delta_M\rangle$ is a condition in $P_S$.
Therefore, $x_{\delta_M}\in S$, and hence $N\cap\kappa\in S$.

We now proceed to prove the claim.
We just need to check that $x_{\delta_M} = M[G]\cap \kappa.$
We have $M\subset M[G]$ and $G\in M[G]$.
In $V[G]$, $G$ defines a bijection $f : \omega_1\to\kappa$.
Let $\dot f\in M$ be a canonical name for this $f$.
We then have that
$$\force \forall p\in \dot
G\;\exists\;\alpha<\omega_1\;\forall\;\gamma<\dom(p)\;p(\gamma)\subset
\dot f''\alpha$$ and
$$\force\forall\;\alpha<\omega_1\;\exists\;p\in\dot
G\;\forall\gamma<\alpha\;\exists\;\beta<\dom(p)\;\dot
f(\gamma)\in p(\beta).$$
Also, $\dot f/G\in M[G]$ and
$\dot f/G\cap M[G] : \delta_M\to M[G]\cap\kappa$
is a bijection.

First we check that $M[G]\cap\kappa\subset x_{\delta_M}.$

Let $\alpha\in M[G]\cap\kappa$. Let $\dot\alpha\in M$ be a name such
that $\forces \dot\alpha <\kappa$ and $\alpha = \dot\alpha/G$.
Then $$\forces \exists p\in \dot G (\dot\alpha\in \bigcup p).$$
Hence $M\models \exists\;\xi<\omega_1\;\dot\alpha\in \dot x_{\xi}.$ Let
$\dot\xi\in M$ be a name for a countable ordinal such that
$$\forces \dot\alpha\in \dot x_{\dot\xi}.$$
Since the semigeneric condition $q$ is in $G$, let $\xi<\delta_M$ be
such that $q\forces \dot\alpha \in \dot x_\xi$. It follows that
$$\alpha = \dot\alpha/G \in (\dot x/G)_\xi\subset x_{\delta_M}.$$

Secondly, we check that $x_{\delta_M}\subset M[G]\cap\kappa.$

Let $\alpha < \delta_M$. Let $\beta\in x_\alpha$. We show that $\beta\in
M[G]$.

Let $p\in M[G]\cap G$ be such that $x_\alpha = p(\alpha)$.
Let $\dot p \in M$ be such that $\dot p/G = p$. Let $\dot \alpha\in M$
be such that $\dot \alpha/G = \alpha$. Let $\dot \xi\in M$ be such that
$$q\force \dot p(\dot \alpha)\subset \dot f''\dot \xi.$$
It follows that $\beta\in M[G]\cap\kappa.$

\end{proof}

As a corollary, if stationary--set--preserving $=$ semiproper, then
projective stationary $=$ spanning. We shall prove the converse later in
this section.

It follows that WRP implies that projective stationary $=$ spanning.
More precisely,

\begin{corollary}
If every local club in $[H_{(2^\kappa)^+}]^\omega$ contains a club, then
every projective stationary set in $\kkk$ is spanning.
\end{corollary}

Looking at the proof of (b), we observe that the club $D$ in the
definition of spanning is  the club that witnesses semiproperness of
$P_S$. If we replace ``club'' by ``local club'', the proof goes through
as before and we get the following characterization of projective
stationary sets.

\begin{lemma}
A set $S\subset\kkk$ is projective stationary if and only if for
every $\lambda\geq \kappa$, for every club $C\subset [\lambda]^\omega$,
there exists a local club $D$ in $[\lambda]^\omega$ such that every
$x\in
D$ has an $\omega_1$--extension $y$ in $C$ such that $y\cap\kappa\in S$.
\end{lemma}

The quantifier $\forall\; C$ in Definition 2.6 and Lemma 4.6 can
be removed by the following trick. Let $S$ be a stationary set in $\kkk$
and let $\lambda\geq \kappa^+$ and $\mu = \lambda^+$. Let

\vskip10truept
\noindent$(4.3)$ \hskip1truecm $S^*_\lambda = \{ M\cap
H_\lambda\;|\;M\in [H_\mu]^\omega,\; S\in
M\;\hbox{\rm and}\;M\cap\kappa\in S\}$, \vskip10truept

\noindent and

\vskip10truept
\noindent$(4.4)$ \hskip1truecm Sub($S^*_\lambda$) $=\{ M\in \hhh\;|\;
M$ has an $\omega_1$--extension $N\in S^*_\lambda\}$.
\vskip10truept

Here we assume that $H_\mu$ has Skolem functions and $M\in
[H_\mu]^\omega$ is an elementary submodel. The set $S^*_\lambda$ is a
stationary subset of $\hhh$ and is equivalent to the lifting of $S$.

\begin{lemma}
(a) $S$ is spanning if and only if Sub($S^*_\lambda$) contains a
club.

(b) $S$ is projective stationary if and only if Sub($S^*_\lambda$)
is a local club.
\end{lemma}

\begin{proof}
We prove (a) as (b) is proved similarly.

Let $\lambda\geq\kappa^+$ and $\mu=\lambda^+$.

First assume that $S$ is spanning. Let
$$C = \{ M\cap H_\lambda\;|\;M\in [H_\mu]^\omega\; \hbox{\rm
and}\;S\in M\}.$$ Let $D$ be a club in $\hhh$ such that every
$M\in D$ has an $\omega_1$--extension $N\in C$ with
$N\cap\kappa\in S$. Then $D\subset $ Sub($S^*_\lambda$).

Conversely, assume that $S$ is not spanning. Let $C = C_F$ be the
least counterexample. As $F$ is definable in $H_\mu$ from $S$, it
belongs to every elementary countable submodel $M$ of $H_\mu$ such
that $S\in M$. Hence every $N\in S^*_\lambda$ is closed under $F$
and it follows that $S^*_\lambda\subset C$. Therefore, every $M\in
$ Sub($S^*_\lambda$) has an $\omega_1$--extension $N\in C$ such
that $N\cap\kappa\in S$. Since $C$ is a counterexample,
Sub($S^*_\lambda$) does not contain a club.
\end{proof}

Now we prove that projective stationary $=$ spanning implies that
stationary--set--preserving $=$ semiproper. This is a consequence of the
following lemma.

\begin{lemma}
Let $P$ be a forcing ($|P|\geq\aleph_1$) and let $\lambda\geq |P|^+$.

(a) $P$ is semiproper if and only if the set $(3.1)$ is spanning.

(b) $P$ preserves stationary sets in $\omega_1$ is and only if the
set $(3.1)$ is projective stationary.
\end{lemma}

\begin{proof}
Both (a) and (b) have the same proof, using Definition 2.6 and Lemma
4.6. The left--to--right implications are obvious, as club implies
spanning and local club implies projective stationary. Thus assume (for
(a)) that the set (3.1) is spanning. If follows from Definition 2.6
that there exists a club $D$ in $\hhh$ such that every $M\in D$ has an
$\omega_1$--extension in the set (3.1). But since every condition that
is semigeneric for an $\omega_1$--extension of $M$ is semigeneric for
$M$, it follows that every $M\in D$ belongs to the set (3.1). Thus the
set (3.1) contains a club and $P$ is semiproper.
\end{proof}

\begin{corollary}
If every projective stationary set is spanning, then every forcing that
preserves stationary sets of $\omega_1$ is semiproper.
\end{corollary}

We conclude Section 4 with the following diagram describing the
implications under the assumption of WRP.

$$\begin{array}{ccccc}
&&{\hbox{\rm full}}&&\cr
&\nearrow&&\searrow&\cr
{\hbox{\rm club $=$ local club}}&&&&{\hbox{\rm projective stationary $=$
spanning}}\cr
&\searrow&&\nearrow&\cr
&&{\hbox{\rm reflective}}&&
\end{array}$$
$$\hbox{\rm Fig.\ } 4.1$$

\section{Strong reflection principle}
The Strong Reflection Principle (SRP) is the statement that every
projective stationary set contains an $\omega_1$--chain. Thus SRP
implies that every projective stationary set is reflective and that
every full set contains a club. As SRP implies WRP (cf. \cite{FJ1}) we
also have local club $=$ club and projective stationary $=$ spanning,
obtaining the diagram (Fig. 1.2) from the introduction.

We shall now look more closely at spanning sets and prove, among others,
that if all spanning sets contain an $\omega_1$--chain then SRP holds.

\begin{definition}
For $X\subset \kkk$, let
$$X^{\bot} = \{ M\in [H_{\kappa^+}]^\omega\;|\;M \hbox{\rm \ has no
$\omega_1$--extension $N$ such that } N\cap \kappa\in X\}.$$
\end{definition}

The set $X^{\bot}$ is a subset of $[H_{\kappa^+}]^\omega$ and is
disjoint from $\hat X$. If $X$ is nonstationary, then $X^\bot$ contains
a club. Let us therefore restrict ourselves to stationary sets
$X\subset \kkk$.

\begin{lemma}
{(i)} If $S_1\subset S_2\subset\kkk$, then $S_2^\bot\subset
S_1^\bot.$

{(ii)} If $S_1\equiv S_2$ mod club filter, then $S_1^\bot\equiv
S_2^\bot$ mod club filter.

{(iii)} $\hat S\cup S^\bot$ is spanning (where $\hat S$ is the
lifting of $S$ to $H_{\kappa^+}$).

{(iv)} $S$ is spanning if and only if $S^\bot$ is nonstationary.
\end{lemma}

\begin{proof}
(ii) Let $F : \kappa^{<\omega}\to \kappa$ be such that $S_1\cap C_F =
S_2\cap C_F$. Let $D = \{ M\in \kkp\;|\;F\in M\}.$ $D$ is a club in
$\kkp$. We claim that $S_1^\bot\cap D = S_2^\bot\cap D.$

If $M\in D$ and $M\not\in S_1^\bot$, then $M$ has an
$\omega_1$--extension $N$ such that $N\cap\kappa\in S_1$. Since $F\in
M\subset N$, $N\cap \kappa$ is closed under $F$. So $N\cap\kappa\in
S_2$. Hence $M\not\in S_2^\bot$. Similarly for the other direction, and
so we have $S_1^\bot\cap D = S_2^\bot\cap D.$

(iii) Let $\lambda\geq \kappa^+$ be arbitrary and let $C$ be a club in
$\hhh$. Let $F : H_\lambda^{<\omega}\to H_\lambda$ be such that
$C_F\subset C$. We claim that every $M\in C_F$ has an
$\omega_1$--extension $N\in C$ such that $N\cap H_{\kappa^+}\in \hat
S\cup S^\bot$, i.e., either $N\cap \kappa\in S$ or $N\cap
H_{\kappa^+}\in S^\bot$.

Let $M\in C_F$. If $M\cap H_{\kappa^+}\in S^\bot$, then we are done.
Otherwise, let $M_0 = M\cap H_{\kappa^+}$. $M_0$ has an
$\omega_1$--extension $N_0\in \kkp$ such that $N_0\cap\kappa\in S$. Let
$N$ be the closure of $M\cup (N_0\cap\kappa)$ under $F$. We have that
$N\in C$ and $M\subset N$. By an argument exactly as in the proof of
Lemma 4.2, we conclude that $N\cap\kappa = N_0\cap\kappa.$ Hence $N$ is
an $\omega_1$--extension of $M$ and $N\cap\kappa\in S$.

(iv) If $S$ is spanning then by definition the set of all $M\in \kkp$
that do have an $\omega_1$--extension $N$ with $N\cap\kappa\in S$
contains a club, and hence $S^\bot$ is nonstationary. If $S^\bot$ is
nonstationary, then, since $\hat S\cup S^\bot$ is spanning, $\hat S$
must be spanning. Hence $S$ is spanning.
\end{proof}

\begin{theorem}
If every spanning set in $\kkp$ contains an $\omega_1$--chain, they
every projective stationary set in $\kkk$ contains an $\omega_1$--chain.
\end{theorem}

\begin{proof}
Let $S$ be a projective stationary set in $\kkk$. By Lemma 5.2(iii),
$\hat S\cup S^\bot$ is spanning in $\kkp$ and therefore contains an
$\omega_1$--chain $\langle M_\alpha\;|\;\alpha<\omega_1\rangle$.  We
claim that
$\{\alpha<\omega_1\;|\;M_\alpha\cap\kappa\in S\}$ contains a club and
therefore $S$ contains an $\omega_1$--chain.

Suppose not. The set $A=\{\alpha<\omega_1\;|\;M_\alpha\in
S^\bot\;\hbox{\rm and}\;\alpha = \delta_{M_\alpha}\}$ is stationary. Let
$$C = \{N\in\kkp\;|\;\kappa\in N\;\hbox{\rm and}\;(\forall\;\beta\in
N\cap\omega_1)\;M_\beta\in N\}.$$
$C$ is a club in $\kkp$. Since $S$ is projective stationary, there
exists an $N\in C$ such that $\delta_N\in A$ and $N\cap\kappa\in S$.
For every $\alpha <\delta_N$ we have $M_\alpha\subset N$. Hence
$M_{\delta_N}\subset N$ and $M_{\delta_N}\cap\omega_1 =
N\cap\omega_1=\delta_N.$ Therefore, $M_{\delta_N}\not\in S^\bot$. This
is a contradiction.
\end{proof}

\begin{corollary}
If every spanning set contains an $\omega_1$--chain, then SRP holds.
\end{corollary}

\section{A structure theorem}
The following definition relativizes projective stationary and spanning.

\begin{definition}
Let $A$ be a stationary set of countable ordinals and let $S\subset
\kkk$.

{(a)} $S$ is {\sl projective stationary above} $A$ if for every
stationary $B\subset A$, the set $\{x\in S\;|\;\delta_x\in B\}$ is
stationary.

{(b)} $S$ is {\sl spanning above} $A$ if for every club
$C\subset\kkp$ there exists a club $D$ in $\kkp$ such that every $M\in
D$ with $\delta_M\in A$ has an $\omega_1$--extension $N\in C$ such that
$N\cap\kappa\in S$.
\end{definition}

The following result is proved in \cite{FJ2}.

\begin{lemma}
If the nonstationary ideal on $\omega_1$ is saturated, then for every
stationary set $S\subset\kkk$ there exists a stationary
$A\subset\omega_1$ such that $S$ is projective above $A$.
\end{lemma}

Notice that the conclusion of the lemma can be stated as: the complement
of $S$ is not full. Thus Lemma 6.2 is a reformulation of Theorem 3.8(a).

\begin{corollary}
If the nonstationary ideal on $\omega_1$ is saturated then  for every
stationary $S\subset\kkk$ there exists a stationary $A\subset
\omega_1$ such that

{(i)} $S$ is projective stationary above $A$, and

{(ii)} $\{x\in S\;|\;\delta_x\not\in A\}$ is nonstationary.
\end{corollary}

\begin{proof}
Let $W$ be a maximal antichain of stationary sets $A\subset \omega_1$
such that $S$ is projective stationary above $A$. Since
$|W|\leq\aleph_1$, there exists a stationary $A_S$ such that
$$A_S = \Sigma\{A\;|\;A\in W\}$$
in the Boolean algebra $P(\omega_1)/{NS}$. It is easy to verify that
$A_S$ has the two properties.
\end{proof}

\begin{corollary}
SRP implies WRP. In fact, assuming SRP, for every stationary
$S\subset\kkk$ there exists a set $X$ of size $\aleph_1$ such that
$\omega_1\subset X$ and an $\omega_1$--chain $\langle
N_\alpha\;|\;\alpha<\omega_1\rangle$ with $\alpha = \delta_{N_\alpha}$
for all $\alpha<\omega_1$ such that $X =
\displaystyle{\bigcup_{\alpha<\omega_1} N_\alpha}$ and $N_\alpha\in S$
for every $\alpha\in A_S$.
\end{corollary}

\begin{proof}
The set $S\cup \{x\;|\;\delta_x\not\in A_S\}$ is projective stationary
and by SRP it contains an $\omega_1$--chain.
\end{proof}

The proof that WRP implies that projective stationary equals spanning
applies to the relativized notions, i.e., projective stationary above
$A$ equals spanning above $A$. Thus we obtain the following theorem.

\begin{theorem}
Assume SRP. Let $\kappa\geq\omega_2$ and let $S\subset\kkk$ be
stationary. There exists a stationary $A_S$ such that

{(i)} for almost all $x\in S,\;\delta_x\in A_S$, and

{(ii)} almost all $x$ with $\delta_x\in A_S$ have an
$\omega_1$--extension $y\in S$.

Moreover, the set $A_S$ is unique mod club filter and if $S_1\equiv
S_2$ then $A_{S_1} \equiv A_{S_2}$.

Also, a stronger version of (ii) holds: for every
$\lambda\geq\kappa$ and every model $(\lambda,\cdots)$, almost all
countable $M\prec (\lambda,\cdots)$ with $\delta_M\in A_S$ have an
$\omega_1$--extension $N\prec(\lambda,\cdots)$ such that
$N\cap\kappa\in S.$
\end{theorem}

\section{Order types and canonical functions}

Two functions $f, g :\omega_1\to\omega_1$ are equivalent (mod club
filter) if the set $\{\alpha<\omega_1\;|\;f(\alpha) = g(\alpha)\}$
contains a club. $f<g$ if and only if
$\{\alpha<\omega_1\;|\;f(\alpha)<g(\alpha)\}$ contains a club. Then $<$
is a well--founded partial order of the equivalence classes and every
function can be assigned a rank in this partial order. For all
$\eta<\omega_2$, there exist {\sl canonical function} $f_\eta$ such that
each $f_\eta$ has rank $\eta$ and when $\eta$ is a limit ordinal then
$f_\eta$ is the least upper bound of $\{f_\xi\;|\;\xi<\eta\}$. The
canonical functions are unique and for $\omega_1\leq \eta<\omega_2$, if
$g_\eta$ is any one--to--one mapping of $\omega_1$ onto $\eta$, then for
almost all $\alpha<\omega_1$,

\vskip10truept
\noindent $(7.1)\hskip1truecm f_\eta(\alpha) = \hbox{\rm order type of }
\{g_\eta(\beta)\;|\;\beta<\alpha\}.$

\vskip10truept
The {\sl Boundedness Principle} is the statement
\vskip10truept

\noindent$(7.2)\hskip1truecm (\forall\;g :
\omega_1\to\omega_1)(\exists\;\eta<\omega_2)\;g<f_\eta.$

\vskip10truept
This follows from the saturation of the nonstationary ideal on
$\omega_1$ (but the consistency strength is considerably less).

\begin{theorem}
The boundedness principle is equivalent to the following statement: for
every club $C\subset\omega_1$, the set
\vskip10truept

\noindent$(7.3)\hskip2truecm \{x\in [\omega_2]^\omega\;|\;\hbox{\rm
order--type}(x)\in C\}$
\noindent is a local club.
\end{theorem}

\begin{proof}
First assume that for every club $C$ the set $(7.3)$ is a local club.
Let $g : \omega_1\to\omega_1$ be an arbitrary function.

Let $C =
\{\gamma<\omega_1\;|\;(\forall\;\alpha<\gamma)\;g(\alpha)<\gamma\}.$

Let $\eta$ and
$\langle x_\alpha\;|\;\alpha<\omega_1\rangle$
be such that $\omega_1<\eta<\omega_2$ and
$\langle x_\alpha\;|\;\alpha<\omega_1\rangle$ is an $\omega_1$--chain
which is a club in $[\eta]^\omega$
and for all $\alpha<\omega_1$
order--type$(x_\alpha)\in C$. By our assumption, such $\eta$ exists.

We claim that $g<f_\eta$. By (7.1), $f_\eta(\alpha) = $
order--type$(x_\alpha)$ for almost all $\alpha<\omega_1$. Let
$$D=\{\alpha\in C\;|\;\alpha <f_\eta(\alpha) = \hbox{\rm
order--type}(x_\alpha)\}.$$
For each $\alpha\in D$ we have $f_\eta(\alpha)\in C$ and
$f_\eta(\alpha)>\alpha$, while $g(\alpha) <\alpha'$, where $\alpha'$ is
the least element of $C$ greater than $\alpha$. Thus $g<f_\eta$,
witnessed by $D$.

Conversely, assume that for every $g : \omega_1\to \omega_1$, there
exists an $\eta<\omega_2$ such that $g<f_\eta$. Let
$C\subset\omega_1$. Consider the set
$$D=\{\eta<\omega_2\;|\;\{\alpha<\omega_1\;|\;f_\eta(\alpha)\in
C\}\hbox{\rm \ contains a club}\}.$$

Using canonicity, it is easy to verify that $D$ is closed. We claim that
$D$ is unbounded.

Let $\eta_0<\omega_2$. We construct a sequence of functions $\langle
g_k\;|\;k<\omega\rangle$ and a sequence of ordinals
$\langle\eta_k\;|\;k<\omega\rangle$ so that
$$f_{\eta_0}<g_0<f_{\eta_1}<g_1<\cdots$$
and that $g_k(\alpha)\in C$ for every $k$ and every $\alpha$. This can
be done since $C$ is unbounded and by our assumption. Let
$$\eta = \hbox{\rm lim}_k\;\eta_k.$$
Then for almost $\alpha$,
$$f_\eta(\alpha) = \hbox{\rm lim}_{k}\; f_{\eta_k}(\alpha)
=\hbox{\rm lim}_{k}\; g_k(\alpha).$$ Since $C$ is closed, we have
$f_\eta(\alpha) \in C $ for almost $\alpha$, and so $\eta\in D.$

Now if $\eta\in D$ and $\langle x_\alpha\;|\;\alpha<\omega_1\rangle$ is
a club in $[\eta]^\omega$, then by
 (7.1) the order type of $x_\alpha$ is $f_\eta(\alpha)$
for almost all $\alpha<\omega_1$, and therefore
$$\{x\in [\eta]^\omega\;|\;\hbox{\rm order--type}(x)\in C\}$$
contains a club in $[\eta]^\omega$. Thus (7.3) is a local club.
\end{proof}

\begin{corollary}
If SRP holds then for every stationary set $S\subset \kkk$,
the set $\{\hbox{\rm order--type}(x\cap\omega_2)\;|\;x\in S\}$ is
stationary.
\end{corollary}

\begin{proof}
SRP implies both the boundedness principle and that local club $=$ club,
and so the set
$$\{x\in \kkk\;|\;\hbox{\rm order--type}(x\cap\omega_2)\in C\}$$
contains a club for every club $C\subset \omega_1.$
\end{proof}

\bibliographystyle{plain}
\bibliography{fjz}

\vskip30truept

\noindent Institute of Mathematics, AMSS, Chinese Academy of Sciences,
Zhong Guan Cun, Beijing 100080, China\\
\noindent{Email:} {\tt qifeng@mail.math.ac.cn}\\
\noindent and \\
\noindent Department of Mathematics,
National University of Singapore, 2 Science Drive 2,
Singapore 117543, Republic of Singapore\\
\noindent{Email:} {\tt matqfeng@math.nus.edu.sg}
\vskip10truept
\noindent Mathematical Institute, The Academy of Sciences of the
Czech
Republic, \v Zitn\'a 25, 115 67 Praha 1, Czech Republic\\
 \noindent{Email:} {\tt jech@math.cas.cz}
\vskip10truept
\noindent Department of Mathematics, University of
Florida, Gainesville,
FL 32611, USA\\
\noindent{Email:} {\tt jinzap@yahoo.com}

\end{document}